# Optimization of Radar Search Patterns for Multiple Scanning Missions in Localized Clutter

Yann Briheche*[†], Frederic Barbaresco*, Fouad Bennis[†], Damien Chablat[†]
* THALES AIR SYSTEMS
Voie Pierre-Gilles de Gennes, 91470 Limours, France
Email: yanis.briheche@thalesgroup.com, frederic.barbaresco@thalesgroup.com
[†]Institut de Recherche en Communications et Cybernetique de Nantes, UMR CNRS 6597
1, rue de la Noe, 44321 Nantes, France
Email: yanis.briheche@irccyn.ec-nantes.fr, fouad.bennis@irccyn.ec-nantes.fr, damien.chablat@irccyn.ec-nantes.fr

*Abstract*—Electronic Phased-Array Radars offer new possibilities for Optimization of Radar Search Pattern by using bi-dimensional beam forming and beam steering, along both elevation and azimuth axes.

The minimization of the Time-Budget required for multiple Radar scanning missions in localized clutter, under constraints of range and detection probability, can be approximated as a Set Cover Problem. We present a Set Cover Problem approximation for Time-Budget minimization of the Radar Search Pattern, and solved this optimization problem using Integer Programming methods based on Branch&Bound and linear relaxation.

## I. INTRODUCTION

Multi-Function Radars (MFR) are a new family of defense systems improving and widening the capabilities of radars [1], [2]. A key feature of MFR is to use the new flexibility offered by electronic scanning to achieve adaptive management of resources [3], [4]. Approximating the Radar Search Pattern optimization problem as a Set Cover Problem offers a flexible yet powerful formulation [5]. Localized clutter and multiple missions can be integrated into this formulation without changing its structure as an optimization problem. This extended formulation can thus be solved with the same methods.

In this paper, we present an extension of the Time-Budget minimization problem to achieve detection for different missions, each with given desired detection range, detection probability, false-alarm probability, Radar Cross-Section (RCS) and Swerling Model.

A few works have been done on the optimization of Radar Scanning coverage using a grid of pencil beams tiled over the entire surveillance space [6], [7]. Those approaches however do not consider the use of different beam-shapes to cover the space. A similar problem is the WiFi Network Cover: for a given base station and given clients, ensure connection for all clients with a minimal numbers of WiFi covers [8], [9]. Radar Search Pattern and WiFi Network Cover have similar underlying structures as optimization problems. Both can be viewed as instances of the Set Cover Problem (SCP). SCP is a classical problem in Combinatorial Optimization, where the objective is to cover a set of elements using a minimum number of available covers. The theoretical problem is known to be generally NP-difficult to solve [10]. In practice, this problem can however be efficiently solved by Branch&Bound exploration using Linear Relaxation for lower bound estimation.

We present the original optimization problem of the Radar Search Pattern to be solved in II-A, and we define the dwell model used to achieve coverage in II-B. We describe our approximation procedure to rewrite the original problem as a SCP in III. Then in III-E, we present an optimization algorithm for the approximated problem. Section IV details implementation and simulation results on a case study and section V concludes on future works and possible improvements of the procedure.

## II. ORIGINAL PROBLEM

Our objective is to find a Radar Search Pattern, i.e. a set of dwells covering and ensuring detection over the surveillance space, while using minimum time-budget.

Surveillance space $\mathcal{A}_S$ is defined as a rectangle in azimuth-elevation space (Fig. 1): $\mathcal{A}_S = [az_{min}, az_{max}] \times [el_{min}, el_{max}]$

### A. Constraints

We must ensure detection for a set of $I$ known missions. For each missions $i \in \{1, \ldots, I\}$, let:
- $\sigma_i$ be the Radar Cross-Section.
- $R_{c,i} : \mathcal{A}_S \to \mathbb{R}^+$ be the desired range of detection.
- $SW_i$ be the Swerling Model.

Let $P_d, P_{fa} \in ]0, 1[$ be the desired detection probability and the desired false alarm probability.

Let $\alpha : \mathcal{A}_S \to [0, 1[$ be the clutter eclipse coefficient. The clutter eclipse coefficient is the ratio of eclipsed area on the range-Doppler map: it corresponds to the probability that a target with random speed at random range is eclipsed in clutter.

Detection is ensured if for each mission $i$ and for each direction $(az, el) \in \mathcal{A}_S$, our Radar Search Pattern contains at least one dwell capable of detecting a target with RCS $\sigma_i$ at range $R_{c,i}(az, el)$ in direction $(az, el)$ with at least detection probability $P_d$ and at most false alarm probability $P_{fa}$ while accounting for clutter eclipse coefficient $\alpha(az, el)$.

### B. Dwell Model

A dwell is the combination of a transmission pattern produced by the radiating elements of the radar, and of a waveform which is the signal sent along the temporal axis.

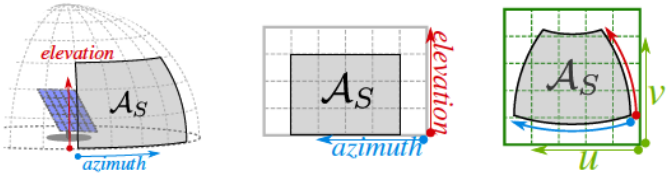

Fig. 1. $\mathcal{A}_S$ in 3D (left), in $az$-$el$ (middle), in direction cosines (right)

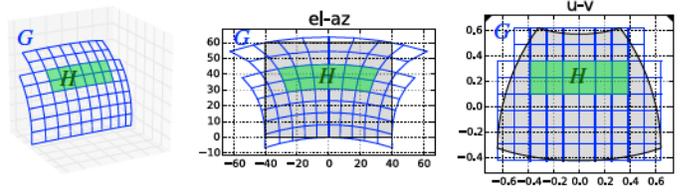

Fig. 2. Grid $G$ and a sub-rectangle $H$ in 3D (left), in $az$-$el$ (middle), in direction cosines (right)

*1) Transmission Pattern:* To achieve detection, we use a 2D linear Phased Array Radar (PAR), containing $K$ rows with spacing $d_y$ and $L$ columns with spacing $d_x$ of isotropic radiating elements. We suppose we can control the phase and amplitude each radiating element:

$$a_{k,l} = A_{k,l}\, e^{j\phi_{k,l}}$$

with $A_{k,l} \in [0,1]$, $\phi_{k,l} \in [0, 2\pi[$, $k \in \{1,...,K\}$ and $l \in \{1,...,L\}$. The theoretical transmission gain of the PAR can be computed from the following formula [11]:

$$g_t(u,v) = \sum_{k=0}^{K-1} \sum_{l=0}^{L-1} a_{k,l}\, e^{j2\pi(kd_y v + l d_x u)/\lambda} \quad (1)$$

with

$$\begin{aligned} u &= \cos(el)\sin(az) \\ v &= \sin(el)\cos(t) - \sin(t)\cos(az)\cos(el) \end{aligned} \quad (2)$$

where $t$ is the PAR tilt angle, $\lambda$ is the signal wavelength, and $(u,v)$ are the direction cosines with values in disk of center $(0,0)$ with radius 1 in $\mathbb{R}^2$.

*2) Waveform:* We consider a given finite set of already designed waveforms $\mathcal{W} = \{w_1, w_2, ..., w_P\}$. For each waveform $w \in \mathcal{W}$, we know the following characteristics:
- time duration $T_w$
- wavelength $\lambda_w$
- SNR detection threshold $s_w(i,\alpha)$ for given:
  - detection and false alarm probability $(P_d, P_{fa})$
  - mission $i$ and associated Swerling Model $SW_i$
  - clutter eclipse coefficient $\alpha$

Informally, we know how much time each waveform requires and how well each waveform can detect a target in clutter for all possible target/clutter situations.

### C. Detection range

The maximum detection range for mission $i$ in direction $(az, el)$ using a dwell $d$ with transmission pattern $g_t$ and waveform $w$ is given by the radar equation [12]:

$$R_{d,i}(az, el)^4 = \frac{P_m\, T_w\, g_t(az,el)\, g_r\, \lambda_w^2\, \sigma_i}{(4\pi)^3\, s_w(i, \alpha(az,el))\, L_u\, L_s(az,el)^2} \quad (3)$$

where $P_m$ is the transmission mean power of the radar, $g_r$ is the reception gain of Digital Beam-Forming, $L_s$ are scanned losses and $L_u$ are uniform losses independent of observation direction. From the radar equation, we can deduce the scanned area of a given dwell $d = (g_t, w)$ for mission $i$ as:

$$\mathcal{A}_{d,i} = \{(az, el) \in \mathcal{A}_S,\ s.t.\ R_{d,i}(az,el) \geq R_{c,i}(az,el)\} \quad (4)$$

For each dwell, the scanned area size in $(u,v)$ is limited by the Digital Beam-Forming bandwidth: $A_d = \iint_{\mathcal{A}_{d,i}} du\, dv \leq A_{\max}$

### D. General formulation

Finding an optimal Radar Search Pattern $\mathcal{S}_{opt}$ is a minimization problem under constraints:

$$\min \sum_{0 \leq j \leq J} T_{w_j} \quad (5a)$$

$$s.t.\ \mathcal{S} = \{d_j = (g_j, w_j), 0 \leq j \leq J\} \quad (5b)$$

$$\forall i \in \{0,...,I\}, \mathcal{A}_S \subset \bigcup_{d \in \mathcal{S}} \mathcal{A}_{d,i} \quad (5c)$$

$$\forall d \in \mathcal{S}, \forall i,\ A_{d,i} = \iint_{\mathcal{A}_{d,i}} du\, dv \leq A_{\max} \quad (5d)$$

Our objective is to find a finite Radar Search Pattern $\mathcal{S}$ (5b), whose dwells cover the entire surveillance area for all missions (5c) and can all be processed at reception (5d), and to minimize its Time-Budget (5a).

## III. SET COVER PROBLEM APPROXIMATION

The problem is difficult to solve under its current form because it has continuous variables (phases and amplitudes of the PAR elements for each dwell) and discrete variables (choice of the waveform for each dwell). The number of variables is proportional to the number of dwells, and is not set. Furthermore, the constraint functions $R_{c,i}$ are generally not convex.

### A. Discrete Grid

To solve this problem, we approximate the surveillance area in $(u,v)$ space using a finite bidimensional grid $G$ (Fig. 2). On this grid, we can test the detection on a cell-by-cell basis with a finite number of cells, instead of working with infinite and continuous set of possible azimuth-elevation directions.

We then divide the problem into two sub-problems: first, the generation of a collection of candidate dwells to cover sub-rectangles of the grid, using Pattern Synthesis; and then we formulate the problem to select a subset of those candidate dwells as a Set Covering Problem.

### B. Pattern Synthesis

Let $H$ be a sub-rectangle of the grid $G$ (Fig. 2). The desired pattern to cover $H$ is defined as

$$g \propto \begin{cases} L_s^2 \max_i \left\{ \dfrac{R_{c,i}^4\, s_w(i,\alpha)}{\sigma_i} \right\} & \text{if } (az,el) \in H \\ 0 & \text{if } (az,el) \notin H \end{cases} \quad (6)$$

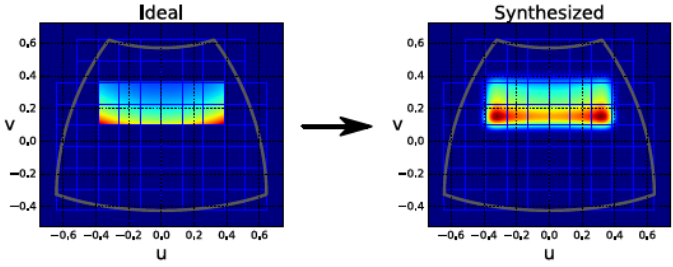

Fig. 3. Ideal pattern (left) and synthesized pattern (right)

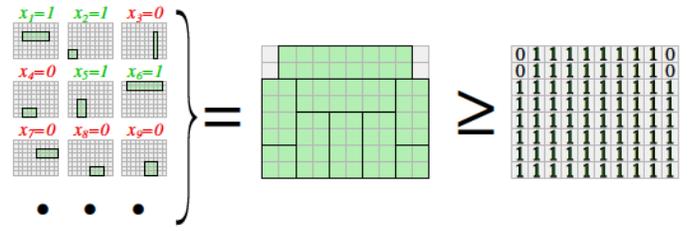

Fig. 5. A set of discrete covers and the sum of its covers for one mission

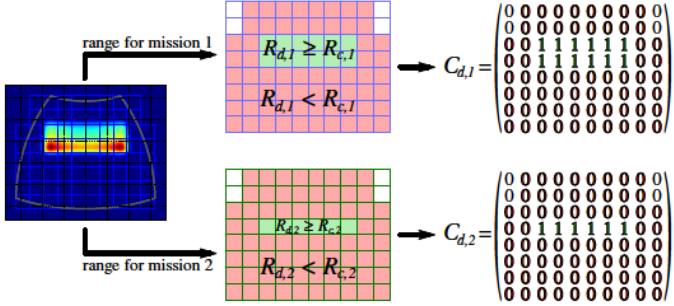

Fig. 4. Discrete covers of a dwell for two scanning missions

by a constant factor, as the PAR feeds are normalized. This pattern provides an ideal distribution for the energetic requirements of the radar surveillance, but is usually unfeasible on a real system.

From the ideal pattern, we synthesize a feasible radiation pattern $g_H$ by applying a 2D Woodward-Lawson sampling method to $g$, adapted from the 1D method presented in [11], [13] (Fig. 3). More sophisticated synthesis methods using least square optimization[14], genetic algorithms[15] and alternating projections[16] are also compatible with our approach.

We compute a collection of radiation patterns, one for each subrectangle on the grid: $\mathcal{T} = \{g_H, H \subset G\}$

### C. Set Cover Problem

The set of candidate dwells $\mathcal{D}$ can be computed as the Cartesian product of $\mathcal{T}$, the set of synthesized patterns, and $\mathcal{W}$, the set of available waveforms:

$$\mathcal{D} = \mathcal{T} \times \mathcal{W} = \{(g_t, w), g_t \in \mathcal{T}, w \in \mathcal{W}\} = \{d_1, \cdots, d_p\}$$

For each dwell $d_j$ in $\mathcal{D}$ and each mission $i$, we compute the discrete cover of the dwell $d_j$ for mission $i$ (Fig. 4):

$$C_{j,i}(m,n) = \begin{cases} 1 & \text{if } \forall (az, el) \in G(m,n), R_{j,i} \geq R_{c,i} \\ 0 & \text{otherwise} \end{cases}$$

where $G(m,n)$ is a cell of $G$; and $R_{j,i}$ is the detection range of dwell $d_j$ for mission $i$. The discrete cover $C_{j,i}$ represents on which cells of the grid the dwell $d_j$ validates the detection range constraint for mission $i$.

Each dell has a different cover for each mission, as energetic requirements differ between missions. Furthermore certain waveforms are more efficient for certain missions.

Then, the problem of finding a Radar Search Pattern ensuring detection over $\mathcal{A}_S$ is equivalent to finding a set of dwells whose sum of discrete covers cover the entire grid $G$ for each mission, such that there is no cell whose value is null (Figure 5).

Since we want to minimize the time budget of the Radar Search Pattern, the cost of each dwell $d_j$ is its associated waveform duration $T_w$. From now on, we will denote the discrete cover cost by $T_j$.

To each dwell $d_j \in \mathcal{D}$, we associate a binary variable $x_j$, indicating if the dwell $d_j$ is in the Search Pattern or not. Then the minimization problem (5) approximated on grid $G$ is equivalent to the following Set Cover Problem:

$$\begin{aligned} \min \quad & \sum_{j=1}^{p} T_j x_j \\ s.t. \quad & \forall i \in \{1, \ldots, I\}, \forall (m,n), \sum_{j=1}^{p} x_j C_{j,i}(m,n) \geq 1 \\ & \forall j \in \{1, \ldots, p\}, x_j \in \{0, 1\} \end{aligned} \quad (7)$$

### D. Integer Program

This Set Cover Problem can be written as Integer Program (IP) by using matrix formulations:
Let $\mathbf{x} = (x_1 \cdots x_p)^T$, let $\mathbf{T} = (T_1 \cdots T_p)^T$ and let

$$\mathbf{A}_i = \begin{pmatrix} C_{1,i}(0,0) & \cdots & C_{p,i}(0,0) \\ C_{1,i}(0,1) & \cdots & C_{p,i}(0,1) \\ \vdots & \ddots & \vdots \\ C_{1,i}(m,n) & \cdots & C_{p,i}(m,n) \\ \vdots & & \vdots \end{pmatrix} \quad \text{and} \quad \mathbf{A} = \begin{pmatrix} \mathbf{A}_1 \\ \vdots \\ \mathbf{A}_i \\ \vdots \\ \mathbf{A}_I \end{pmatrix}$$

Then we can write our Set Cover Problem as the following IP:

$$\begin{aligned} \min \quad & \mathbf{T}^T \cdot \mathbf{x} \\ s.t. \quad & \mathbf{A} \cdot \mathbf{x} \geq \mathbf{1} \\ & \mathbf{x} \in \{0,1\}^p \end{aligned}$$

### E. Optimization Algorithm

The Integer Program obtained for a multiple missions in localized clutter has the same formulation than single mission surveillance problem without clutter [5], and thus the same optimization method can be used:
IP are generally difficult to solve in a straightforward manner. However the linear relaxation of the problem (where the $x_j$ become continuous variables in $[0, 1]$) can be quickly solved. We used Branch&Bound to explore the solution space, with lower bound generation by Linear Relaxation to avoid enumeration of all possible solutions [17].

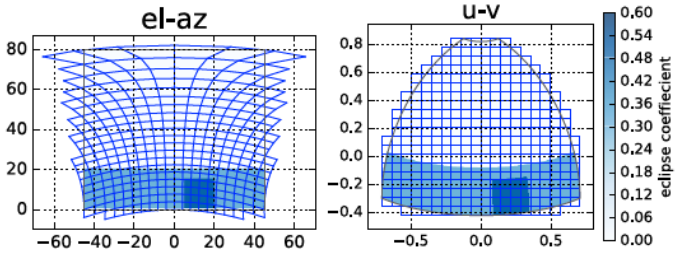

Fig. 6. Eclipse coefficient value on the surveillance area

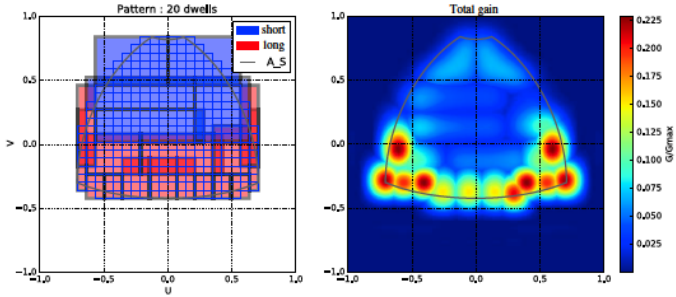

Fig. 7. Optimized solution

## IV. EXPERIMENTAL RESULTS

We applied our procedure for optimization of two scanning missions in presence of localized clutter (Fig. 6). Desired detection range for each mission is delimited by a minimum height and a minimum distance (Fig. 8).

The PAR has $20{\times}20$ half-spaced radiating elements. The grid $G$ is laid on a $20{\times}20$ lattice. We used a set with two possible waveforms $\mathcal{W} = \{w_1, w_2\}$, with a long waveform $w_1$ and a short waveform $w_2$. We computed 32810 feasible dwells in our study case. The grid covering the surveillance area has 360 cases, with 2 scanning missions. The corresponding IP has 32810 variables and 720 inequality constraints.

The computation of the IP is done in Python, and its optimization is done using CPLEX. The total time required to find the solution is 15 seconds on an i7-3770@3.4GHz processor.

The obtained solution uses 20 dwells to cover the surveillance area (Fig. 7), dwells covering the localized clutter have smaller beamwidth, focusing the energy, as they must achieve a higher detection range, and thus require more energy, while dwells at high elevations use the short waveform (in blue).

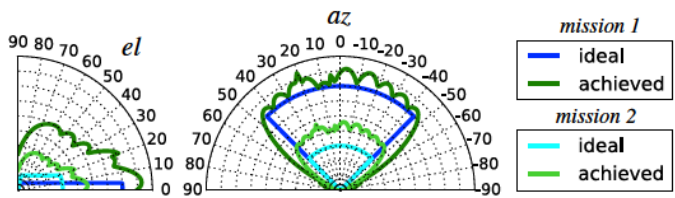

Fig. 8. Detection range of the solution

## V. CONCLUSIONS

The Radar Search Pattern Optimization procedure presented in this paper is an improved formulation which can account for multiple missions and localized clutter. Multiple mission optimization is solved as one global problem, rather than different individual sub-problems, thus optimizing resource allocation for all missions. Furthermore the algorithm is fast enough to be used in operational situation for real-time optimization of resources.

The flexibility of this framework also allows the integration of other types constraints. Our future works will focus on integration of scan update rates constraints and terrain masking.

## ACKNOWLEDGEMENT

This work is partly supported by a DGA-MRIS scholarship.